\documentstyle[12pt]{article}
\begin{document}
\centerline{\bf Defects for Ample Divisors of Abelian Varieties, Schwarz Lemma,}
\centerline{\bf and Hyperbolic Hypersurfaces of Low Degrees}

\bigskip
\centerline{Yum-Tong Siu and Sai-Kee Yeung\ %
\footnote{Both authors were partially supported by grants from the National Science Foundation
and the second author was also partially supported by a grant from
the Sloan Foundation.}
}

\bigskip
The main purpose of this paper is to prove the following theorem on the 
defect relations for ample divisors of abelian varieties.

\medskip
\noindent
{\it Main Theorem.}  Let $A$ be an abelian variety of complex dimension $n$
and $D$ be an ample divisor in $A$.  Let $f:{\bf C}\rightarrow A$ be
a holomorphic map.  Then the defect for the map $f$ and the divisor $D$
is zero.

\medskip
\noindent
{\it Corollary to Main Theorem.} The complement of an ample divisor $D$ in an abelian variety
$A$ is hyperbolic in the sense that there is no nonconstant holomorphic map
from $\bf C$ to $A-D$.

\medskip
We give also in this paper 
the following results in hyperbolicity problems: (i) a 
Schwarz lemma for general jet differentials and (ii) examples of hyperbolic
hypersurfaces of low degree.

\medskip
\noindent
{\it Theorem 1 (General Schwarz Lemma).} Let $X$ be a compact complex
subvariety in ${\bf P}_N$ and $f:{\bf C}\rightarrow X$ be
a nonconstant holomorphic map and $\zeta$ be the coordinate of $\bf C$.
Let $\omega$ be a holomorphic $k$-jet differential of weight $m$ 
on $X$ which vanishes on ample divisor in $X$.  Then the pullback 
of $\omega$ by $f$ (by which is meant only the part 
of the form $\tau(\zeta)d\zeta^m$) vanishes identically on $\bf C$.  
Here a holomorphic $k$-jet differential on the subvariety $X$
means that there are holomorphic jet differentials on open
coordinate charts in ${\bf P}_N$ such that,
for any holomorphic map $g:U\rightarrow X$ from an open subset $U$ to $X$,
their pullbacks by $g$ agree on $U$.

\medskip
\noindent
{\it Theorem 2 (General Schwarz Lemma for Log-Pole Differentials).} 
Let $X$ be a compact complex
subvariety in ${\bf P}_N$.  Let $Z_1,\cdots,Z_p$ be distinct
irreducible complex hypersurfaces
in $X$ and $f:{\bf C}\rightarrow X-\cup_{j=1}^pZ_j$ be
a nonconstant holomorphic map and $\zeta$ be the coordinate of $\bf C$.
Let $\omega$ be a meromorphic $k$-jet differential on $X$ of at most log-pole
singularity along $\cup_{j=1}^pZ_j$ such that the weight of $\omega$ is
$m$ and $\omega$ vanishes on ample divisor in $X$.  Then the pullback 
of $\omega$ by $f$ vanishes identically on $\bf C$.  
Here a a meromorphic $k$-jet differential on $X$ of at most log-pole
singularity along $\cup_{j=1}^pZ_j$ means that locally it is a polynomial
with constant coefficients whose variables are local holomorphic jet
differentials and meromorphic differentials of the form $d^\nu\log g$, where
$\nu$ is a positive integer and $g$ is a local holomorphic function whose
zero-set is contained in $\cup_{j=1}^pZ_j$.

\medskip
\noindent
{\it Theorem 3}.  Let 
$N=4n-3$ and $p=1+N(N-2)=16(n-1)^2$.
Then for generic linear functions $H_j(x_0,\cdots,x_n)$ on ${\bf C}^{n+1}$ 
($1\leq j\leq N$) the hypersurface $\sum_{j=1}^N H_j^p=0$ 
in ${\bf P}_n$ is hyperbolic.

\medskip
\noindent
{\it Theorem 4}.  Let $n\geq 11$ and $g(x_0,x_1,x_2,x_3)$ be a 
homogeneous polynomial of degree $2$ with
$g(0,0,0,x_3)=x_3^2$ satisfying the condition that,
for $h(\xi,\eta)$ equal to 
$g((-1)^{1\over n}\xi,\xi,\eta,1)$,$g(\eta,(-1)^{1\over n}\xi,\xi,1)$,
or $g(\eta,\xi,(-1)^{1\over n}\xi,1)$ with $(-1)^{1\over n}$ equal to any 
$n^{\hbox{th}}$ root of $-1$,
the polynomial
$$
-\eta^n+{1\over 2}\left({\partial h\over\partial\xi}
\right)^2\left({\partial^2h\over\partial\xi^2}\right)^{-1}-h
$$
of degree $n$ in $\eta$ has $n$ distinct roots
and $\partial^2h\over\partial\xi^2$ is a nonzero constant.
Then the surface defined by 
$$
x_0^n+x_1^n+x_2^n+x_3^{n-2}g(x_0,x_1,x_2,x_3)=0
$$
is hyperbolic.

\medskip
\noindent
{\it Corollary to Theorem 4}.  
Let $n\geq 11$ and $a_0,a_1,a_2$ be complex numbers.
Suppose $a_i^n\not=(-1)^{n+1}a_j^n$ for $0\leq i<j\leq 2$
and $1+a_j({-2a_j\over n})^{2\over n-2}+({-2a_j\over n})^{n\over n-2}\not=0$
for $0\leq j\leq 2$ and for all $(n-2)^{\hbox{th}}$ roots in that condition.
Then the surface defined by 
$$
x_0^n+x_1^n+x_2^n+x_3^{n-2}(x_3^2+a_0x_0^2+a_1x_1^2+a_2x_2^2)=0
$$
is hyperbolic.

\medskip
The Corollary to the Main Theorem was already proved in [SY96b].  The proof in [SY96b] could 
yield the result that the defect for $f$ and $D$ is less than $1$, but the 
method there could not yield the sharp result of the defect being zero.

\medskip
The idea in the proof of the Main Theorem is as follows.  First of all
we can assume without loss of generality that the image of $f$ is Zariski
dense in $A$.  For any given positive number
$\epsilon$ we construct a meromorphic $k$-jet differential $\omega$
on $A$ whose pole is dominated by a pole divisor of order $\delta$ along $D$ 
such that for some positive integer $q$ with ${\delta\over q}<\epsilon$ 
the $k$-jet differential 
$\omega$ vanishes to order at least $q$ along the $k$-jet
space $J_k(D)$ of $D$ when $\omega$ is regarded as a function on the
$k$-jet space $J_k(A)$ of $A$. The existence of $\omega$ follows from the
theorem of Riemann-Roch and the fact that the codimension of $J_k(D)$ in
$J_k(A)$ increases without bounds as $k$ increases without bounds. 
We then pull back $\omega$ by $f$ and apply the
logarithmic derivative lemma to the meromorphic function on $\bf C$
defined by the pullback.  In this process we have to use the 
result in [SY95] on the translational invariance of the Zariski closure of
the image of the $k$-order differential of $f$ to make sure that we can
construct an $\omega$ whose pullback by $f$ is not identically zero.

\medskip
Since Ahlfors [A41] introduced his more geometric view to study value
distribution theory, the Schwarz lemma has been one of the indispensable
tools in value distribution theory.  Green and Griffiths [GG79] introduced
a general Schwarz lemma to give an alternative approach to Bloch's theorem
[B26] and gave a sketch of the proof of their general Schwarz lemma.
(For a discussion of the various ideas and proofs for Bloch's theorem see
for example [S95].)
After some unsuccessful attempts to give the details of the Schwarz lemma 
given there, some authors [SY96a, D95, DL96] introduced 
jet differentials with special properties to get a proof for the Schwarz
lemma.  Because of such unsuccessful attempts, there were some skepticisms
(now known to be unjustified) as to the completeness of the version of
the proof of Bloch's theorem presented in [GG79].  Since about two years 
ago, the authors and a number of other people working on hyperbolicity 
problems began to believe
that the Schwarz lemma for general jet differentials 
either exactly or essentially as stated 
in [GG79] could be proved.  Theorem 1 is 
slightly more general than the Schwarz lemma for general jet differentials 
stated in [GG79].  Our proof is completely different from the sketch of the 
proof given in [GG79].  We choose
to present Theorem 1 together with the Main Theorem in this paper, because
the Main Theorem is closely linked to Bloch's theorem [B26] and the general
Schwarz lemma was introduced in the paper of Green and Griffiths [GG79]
which used it for an alternative approach to Bloch's theorem [B26].  
The main idea of the proof of Theorem 1 given here is that 
Cauchy's integral formula for derivatives gives domination of the higher
order derivatives at a point by the first-order derivative on a circle 
centered at that point.  Technically, to handle the complications arising
from the radius of the circle and from different coordinate charts of 
the target manifold, one uses the techniques of value distribution
theory and a simple curvature argument introduced in [S87].  The proof of
Theorem 1 after some very minor modifications works also for meromorphic
jet differentials with only log-pole singularities when the image of the
map from $\bf C$ is disjoint from the log-pole singularities.  Theorem 2
is the log-pole case of Theorem 1.

\medskip
For hyperbolicity problems, the three approaches of the Borel lemma, the
Schwarz lemma for jet differentials linear in the highest order,
and meromorphic connections are very closely related.  The approach of
the Borel lemma, when applicable, usually gives sharper and cleaner arguments
than the other two.  
However, for jet differentials
with higher degree in the highest order, there is a possibility of getting
sharper results but the difficulty to conclude algebraic dependency from
differential equations of the vanishing the pullback of the jet differential
is the main problem which in most cases is insurmontable.  In the final part
of this paper we construct hyperbolic hypersurfaces
of low degrees in any dimension by using the Borel lemma and a simple 
dimension counting argument for certain subvarieties in the Grassmannians.
A variation of the construction can make the degree lower
in the case of a surface.  To prepare for the variation of the construction
we introduce a generalized Borel lemma (Proposition 2)
and, to illustrate the close
relation between the approach of the Borel lemma and that of the jet
differentials, we use the Schwarz lemma to prove the generalized Borel lemma.

\medskip
Examples of hyperbolic hypersurfaces were constructed by Brody-Green [BG77], 
Zaidenberg [Z89], Nadel [89] in dimension 2 and by Adachi-Suzuki [AS90] in
some low dimensions, and finally by Masuda and Noguchi [MN94] in any
dimension.  The degree of a hyperbolic hypersurface constructed by
them is exceedingly high relative to its dimension and the algorithm for the
construction is rather involved.  The degree of the hyperbolic hypersurface
constructed in Theorem 3 is only of the order of the square of
its dimension.  Recently El Goul [E96] gave a construction of a hyperbolic surface
of degree $14$ and it was brought to our attention that a suggestion by 
Demailly could lower the degree in El Goul's construction to $11$.  In 
Theorem 4 we point out how El Goul's construction fits as a variation
in the framework of using the Borel approach to get examples of hyperbolic 
hypersurfaces.  The surface of degree $11$ 
constructed in Theorem 4 differs only very
slightly from the construction of El Goul and Demailly, but we use the approach
of the Borel lemma which, though very closely related to the approach of
meromorphic connections, is simpler, more powerful, and more elegant. 

\bigskip
\noindent
{\bf \S1. Defects for ample divisors of abelian varieties}.

\medskip
To prepare for the proof of the Main Theorem, we now introduce 
some notations and terminology.
For a complex manifold $X$ we denote by $J_k(X)$ the space of 
all $k$-jets in $X$ so that every element of $J_k(X)$ is represented by
$\left({d^\alpha\over{d\zeta^\alpha}}g\right)(0)$ ($0\leq\alpha\leq k$)
for some holomorphic map $g$ from an open 
neighborhood of the origin in $\bf C$ (with coordinate $\zeta$) to $X$.
In particular $J_1(X)$ means the tangent bundle of $X$.
By a holomorphic (meromorphic) $k$-jet differential $\omega$ of 
weight $m$ on an open subset of $G$ of $X$ with local coordinates
$z_1,\cdots,z_n$ we mean an expression of the form
$$
\omega = \sum_{\tilde\nu}\omega_{\nu_{1,1}\cdots\nu_{1,k}\cdots\nu_{n,1}\cdots\nu_{n,k}}
(dz_1)^{\nu_{1,1}}\cdots(d^kz_1)^{\nu_{1,k}}\cdots 
(dz_n)^{\nu_{n,1}}\cdots(d^kz_n)^{\nu_{n,k}}\leqno{(1)}
$$
where the summation is over the $kn$-tuple
$$
\tilde\nu=(\nu_{1,1},\nu_{1,2},\cdots,\nu_{1,k},\cdots,\nu_{n,1},\nu_{n,2},\cdots,\nu_{n,k}) 
$$
of nonnegative integers with
$$
(\nu_{1,1}+2\nu_{1,2}+\cdots+k\nu_{1,k})+\cdots+(\nu_{n,1}+2\nu_{n,2}+\cdots+k\nu_{n,k})=m
$$
and $\omega_{\nu_{1,1}\cdots\nu_{1,k}\cdots\nu_{n,1}\cdots\nu_{n,k}}$ is 
a holomorphic (meromorphic) function on $G$.
For a holomorphic map $g$ from an open subset $U$ of $\bf C$ with
coordinate $\zeta$ to $X$ and for 
any meromorphic $k$-jet differential $\omega$ of weight $m$ on $X$, 
by $g^*\omega$ we mean only the 
term containing $(d\zeta)^m$ or, when there is no confusion, we mean only
the coefficient of the term containing $(d\zeta)^m$. According to this
convention $g^*\omega$ can be regarded as a function on $U$.

\medskip
Because of Bloch's theorem [B26, GG79, K80, McQ96, NO90, Oc77] which 
states that the Zariski closure of the
image of a holomorphic map from $\bf C$ to an abelian variety must be
equal to the translate of an abelian subvariety, to prove the Main Theorem
we can assume without loss of generality that the image of $f$ is Zariski
dense in $A$.  
We denote by $d^kf$ the map from $\bf C$ to $J_k(A)$ induced by $f$.
Denote by $\overline{J_k(A)}$ 
the compactification $A\times{\bf P}_{nk}$ of $J_k(A)=A\times{\bf C}^{nk}$.
Let $H_{kn}$ be the pullback to ${\overline{J_k(A)}}$
of the hyperplane section line bundle of ${\bf P}_{nk}$ by the 
projection map 
$\overline{J_k(A)}=A\times{\bf P}_{nk}\rightarrow{\bf P}_{nk}$.
Let $\pi:{\overline{J_k(A)}}\rightarrow A$ be
the projection onto the base manifold. 
By the Zariski closure of ${\hbox{Im }}d^kf$ in $J_k(A)$ we mean the
intersection with $J_k(A)$ of the Zariski closure of ${\hbox{Im }}d^kf$ in
$\overline{J_k(A)}$.
By [SY96b] the Zariski closure in $J_k(A)$ of the image of 
$d^kf$ is translational invariant
and is therefore of the form  $A\times W_k$ for some irreducible subvariety
$W_k$ of positive dimension in ${\bf P}_{nk}$.

\medskip
Write $A={\bf C}^n/\Lambda$.
Let $L_D$ be the line bundle over
$A$ associated to $D$.  
Since $D$ is ample in $A$, for any integer $p\geq 2$ 
the global holomorphic sections of the line 
bundle $L_D^{\otimes p}$ over $A$ generates the $(p-2)$-jets of $A$.

\medskip
Let $\theta_D$ be the theta function on the universal cover 
${\bf C}^n$ of $A$ which defines the divisor $D$.  We denote by
$J_k(D)$ the subvariety of $J_k(A)$ defined by 
$d^j\theta_D=0$ for $0\leq j\leq k$.  Note that when $D$ is nonsingular,
this notation $J_k(D)$ agrees with the earlier definition of $J_k(X)$ with
$X=D$.  To prove the Main Theorem we can assume without loss of generality
that $L_D=(L^\prime_D)^{\otimes p}$
for some integer $p\geq k+2$, because we can simply replace the lattice
$\Lambda$ defining $A$ by $p\Lambda$.  Let $\Delta$ denote the open unit
disk in $\bf C$ centered at the origin.  Let $k\geq n$.

\medskip
\noindent
{\it Lemma 1.} There exists a holomorphic deformation $D(t)$ ($t\in\Delta$) 
of $D$ such that for $t\in\Delta-0$
the subvariety $J_k(D(t))\cap(A\times W_k)$ is of codimension 
at least $n+1$ in $A\times W_k$.

\medskip
\noindent
{\it Proof.}  Let $V=\Gamma(A,L_D)^{\oplus 2}$.  Since 
$\Gamma(A,L_D)$ generates the $k$-jets of $A$, for every
point $P\in A$ there exist an open neigbhorhood $U_P$ of $P$ in
$A$ and an element $v_P=(v_{P,0},v_{P,1})\in V$ such that

\medskip
\noindent
(i) $[v_{P,0},v_{P,1}]|U_P$ defines a holomorphic map $\Phi_P$
from $U_P$ to the projective line ${\bf P}_1$, and

\medskip
\noindent
(ii) $J_k(\Phi_P^{-1}(Q))\cap(U_P\times W_k)$ is a subvariety of codimension at
least $n+1$ in $U_P\times W_k$ for $Q$ in the projective line ${\bf P}_1$.

\medskip
We can choose a compact neighborhood $K_P$ of $P$ in
$U_P$ such that for a finite number of points $P_1,\cdots,P_\ell$ in
$A$ we have $A=\cup_{j=1}^\ell K_{P_j}$.  Let $E_j$ be the subset
of $V$ consisting of all $v=(v_0,v_1)\in V$ such that the following 
two conditions do not simultaneously hold.

\medskip
\noindent
(i) $(v_0,v_1)$ is nowhere zero on $K_{P_j}$.

\medskip
\noindent
(ii) the codimension of $J_k(\Phi^{-1}(Q))\cap((A-\{v=0\})\times W_k)$ is a 
subvariety of codimension at least $n+1$ in $(A-\{v=0\})\times W_k$ for 
at every point of $K_{P_j}$ for every $Q$ in ${\bf P}_1$, 
where $\Phi:A-\{v\not=0\}\rightarrow{\bf P}_1$ is
defined by $[v_0,v_1]$.

\medskip
Then $E_j$ is a nowhere dense closed subset of $V$.  By Baire category 
theorem we conclude that there exists $v=(v_0,v_1)\in V-\cup_{j=1}^\ell E_j$.
In particular, we conclude that $J_k(\{v_0=0\})\cap(A\times W_k)$ is of 
codimension at least $n+1$ at every point of $A\times W_k$.  We need only
now consider the deformation given by $\theta_{D_t}=\theta_D+tv_0$.  Then
there exists a positive number $\eta$ such that
$J_k(\theta_{D_t})\cap(A\times W_k)$ is of 
codimension at least $n+1$ at every point of $A\times W_k$
for $0<|t|<\eta$.  Q.E.D.

\medskip
Take a fixed small positive integer $\delta$ which can actually be chosen
to be $1$.  We keep the symbol $\delta$ to show the role played by it.
Take a positive integer $q$ which will be very large compared to $\delta$
and then take a positive integer $m$ which will be very large compared to
$q$.  The conditions on the sizes of $q$ and $m$ will be specified later.
We are going to construct a non identically zero
$L_D^{\otimes\delta}$-valued
holomorphic $k$-jet differential of weight $m$ which vanishes to order
at least $q$ along $J_k(D)$ for some sufficiently large $m$.  We first use
the theorem of Riemann-Roch to do this when $D$ is replaced by $D_t$ 
for $t\in\Delta-0$ close to $0$ with $D_t$ satisfying some additional
tranversality condition.  Then we use the semicontinuity of the
dimension of the space of holomorphic sections of line bundles in a 
holomorphic deformation to get the conclusion for $D$ when $t\rightarrow 0$. 
In the following lemma, for notational simplicity
the dependence of the constants on $A$ and $D$ is
not explicitly stated out.

\medskip
\noindent
{\it Lemma 2.}  There exists a postive integer $m_0(W_k,\delta,q)$ 
depending
on $W_k,\delta,q$ (and $A$ and $D$) 
such that for $m\geq m_0(W_k,\delta,q)$ there exists 
an $L_D^{\otimes\delta}$-valued holomorphic
$k$-jet differential on $A$ of weight $m$ whose restriction to
$A\times W_k$ is not identically zero and which vanishes along 
$J_k(D)\cap(A\times W_k)$ to order at least $q$.  In particular, 
from the definition of $W_k$ one knows that $\omega$
is not identically zero on $d^kf$.

\medskip
\noindent
{\it Proof.} 
From Lemma 1 there exists a holomorphic family of ample 
divisors $D(t)$ ($t\in\Delta$) with $D(0)=D$
such that for $t\in\Delta-0$
the subvariety $J_k(D(t))\cap(A\times W_k)$ is of codimension 
at least $n+1$ in $A\times W_k$.
From the ampleness of $D(t)$ and Kodaira's vanishing theorem
we have a positive number $m^\prime_0(W_k,t)$ such that
$$
H^p\left(A\times W_k,\pi^*(L_{D(t)}^{\otimes\delta})\otimes H_{kn}^{\otimes m}|(A\times W_k)\right)=0
$$
for any positive integer $p$ and for $m\geq m^\prime_0(W_k)$.
Let $d$ be the complex dimension of $W_k$.
By the theorem of Riemann-Roch and the K\"unneth formula, we have
$$
\dim_{\bf C}\Gamma\left(A\times W_k,\pi^*(L_{D(t)}^{\otimes\delta})\otimes H_{kn}^{\otimes m}|(A\times W_k)\right)
\geq C_{W_k,t}\delta^n m^d
$$
where $C_{W_k,t}$ is a positive constant dependent on $W_k$
but independent of $m$.
We now choose some $t\in\Delta-0$.
Let ${\cal O}$ be the structure sheaf of $A\times W_k$ and
${\cal I}(t)$ be the ideal sheaf of $\overline{J_k(D(t))}\cap(A\times W_k)$.  
Since for $t\in\Delta-0$ the subvariety $\overline{J_k(D(t))}\cap(A\times W_k)$ of $A\times W_k$ 
is of complex dimension at most $d-1=(n+d)-(n+1)$,  
we have the following estimate
$$
\dim_{\bf C}\Gamma\left(A\times W_k,{\cal O}/{\cal I}(t)^{q+1}\otimes
\left(\pi^*(L_{D(t)}^{\otimes\delta})\otimes H_{kn}^{\otimes m}\right)
\right)
\leq C^\prime_{W_k,\delta,q,t}m^{d-1}
$$
where $C^\prime_{W_k,\delta,q,t}$ is a positive constant dependent on
$W_k,\delta,q,t$ but independent of $m$.  
Thus there exists a positive integer
$m^{\prime\prime}_0(W_k,\delta,q,t)$ such that for 
$m\geq m^{\prime\prime}_0(W_k,\delta,q,t)$ the dimension of
$$
\Gamma\left(A\times W_k,{\cal I}(t)^{q+1}\otimes
\left(\pi^*(L_{D(t)}^{\otimes\delta})\otimes H_{kn}^{\otimes m}\right)
\right)
$$
nonzero.
Clearly we can assume that $m^{\prime\prime}_0(W_k,\delta,q,t)$
is lower semicontinuous as a function of $t\in\Delta-0$.
Let ${\cal F}_{\delta,q,m}$ be the torsion-free sheaf of rank $1$ 
over $A\times W_k\times\Delta$
whose restriction to $A\times W_k\times t$ for a generic $t\in\Delta$ is 
equal to 
${\cal I}(t)^{q+1}\otimes
\left(\pi^*(L_{D(t)}^{\otimes\delta})\otimes H_{kn}^{\otimes m}\right)$.
Let ${\cal G}_{\delta,q,m}$ be locally free sheaf over $\Delta$ which is
the zeroth direct image of the sheaf ${\cal F}_{\delta,q,m}$ under the 
map $A\times W_k\times\Delta\rightarrow\Delta$ which is the projection
onto the last factor.
Fix $t_0\in\Delta-0$. Then for 
$m\geq m^{\prime\prime}_0(W_k,\delta,q,t_0)$ the rank of 
${\cal G}_{\delta,q,m}$ is positive.  Thus for a generic $t\in\Delta$
the dimension of 
$$ 
\Gamma\left(A\times W_k,{\cal I}(t)^{q+1}\otimes
\left(\pi^*(L_{D(t)}^{\otimes\delta})\otimes H_{kn}^{\otimes m}\right)
\right)
$$
is nonzero.  By the semicontinuity of
$$
\dim_{\bf C}
\Gamma\left(A\times W_k,{\cal I}(t)^{q+1}\otimes
\left(\pi^*(L_{D(t)}^{\otimes\delta})\otimes H_{kn}^{\otimes m}\right)
\right)
$$
as a function of $t$ and by letting $t\rightarrow 0$, we conclude 
that for 
$$
m\geq m^{\prime\prime}_0(W_k,\delta,q,t_0)
$$
there exists a non identically zero global holomorphic section
$\omega^\prime$ of 
$\pi^*(L_D^{\otimes\delta})\otimes H_{kn}^{\otimes m}|(A\times W_k)$ 
over $A\times W_k$ which vanishes along $J_k(D)\cap(A\times W_k)$ to 
order at least $q$.
  
\medskip
Let ${\cal I}_{W_k}$ be the ideal sheaf of $W_k$ in ${\bf P}_{nk}$. 
By K\"unneth's formula, from the ampleness of $D$ we have a positive integer
$m^{\prime\prime\prime}_0(W_k)$ such that
$$
H^1\left({\overline{J_k(A)}},\pi^*(L_D^{\otimes\delta})\otimes H_{kn}^{\otimes m}
\otimes{\cal I}_{W_k}\right)=0
$$
for $m\geq m_0^{\prime\prime\prime}(W_k)$.
Thus for $m\geq m_0^{\prime\prime\prime}(W_k)$
we can extend $\omega^\prime$ to an
$L_D^\delta$-valued holomorphic $k$-jet differential $\omega$ on $A$.
Now we need only set $m_0(W_k,\delta,q)$ to be at least as large
as $m_0^\prime(W_k)$, $m^{\prime\prime}(W_k,\delta,q,t_0)$, and
$m^{\prime\prime\prime}_0(W_k)$. Q.E.D.

\medskip
\noindent
{\it Proof of the Main Theorem.} Let $\tilde\omega=\theta_D^{-\delta}\omega$.  
Then $\tilde\omega$ is a 
meromorphic $k$-jet differential on $A$. 
We pull back the meromorphic $k$-jet differential $\tilde\omega$ by $f$ and
recall that by our convention 
$f^*\tilde\omega$ simply means the coefficient of $(d\zeta)^m$ where
$m$ is the weight of $\tilde\omega$ and $\zeta$ is the global coordinate of
$\bf C$.  In other words, we consider $\tilde\omega$
as a function on the space $J_k(A)$ of $k$-jets and $f^*\tilde\omega$ as the evaluation
of $\tilde\omega$ on the image of $d^kf$ when $d^kf$ is considered as a map
from ${\bf C}$ to $J_k(A)$.  
In this sense $f^*\tilde\omega$ is a meromorphic
function on ${\bf C}$.  Write $\tilde\omega$ in terms of the global coordinate 
system and we get
$$
\tilde\omega=\sum_{(\ell_1,\nu_1),\cdots,(\ell_n,\nu_n)}{\tilde a}_{(\ell_1,\nu_1),\cdots,(\ell_n,\nu_n)}
(z_1,\dots,z_n)(d^{\ell_1}z_1)^{\nu_1}\cdots(d^{\ell_n}z_n)^{\nu_n}\leqno{(1)}
$$
where 
${\tilde a}_{(\ell_1,\nu_1),\cdots,(\ell_n,\nu_n)}(z_1,\dots,z_n)$ is a meromorphic function on $A$
with the property that 
$$
\theta_D^\delta{\tilde a}_{(\ell_1,\nu_1),\cdots,(\ell_n,\nu_n)}
$$ is an entire function on ${\bf C}^n$ when ${\tilde a}_{(\ell_1,\nu_1),\cdots,(\ell_n,\nu_n)}$
is used to denote also its pullback to the universal cover ${\bf C}^n$ of $A$.
Let 
$$
a_{(\ell_1,\nu_1),\cdots,(\ell_n,\nu_n)}(z_1,\dots,z_n)
=\theta_D^\delta{\tilde a}_{(\ell_1,\nu_1),\cdots,(\ell_n,\nu_n)}(z_1,\dots,z_n).
$$
Let $\varphi$ be a nonnegative quadratic form on ${\bf C}^n$ such that
$|\theta_D|^2\exp(-\varphi)$ is a well-defined scalar function on $A$.
In other words, $\exp(-\varphi)$ defines a Hermitian metric on the fibers
of the line bundle $L_D$.  Thus 
$$
\exp\left(-{\delta\over 2}\varphi\right)|a_{(\ell_1,\nu_1),\cdots,(\ell_n,\nu_n)}(z_1,\dots,z_n)|
$$
is a well-defined smooth function on $A$ and is therefore bounded.  Hence
$$
|a_{(\ell_1,\nu_1),\cdots,(\ell_n,\nu_n)}(z_1,\dots,z_n)|\leq C
\exp\left({\delta\over 2}\varphi\right)\leqno{(2)}
$$
on ${\bf C}^n$ for some constant $C$.  From (1) we have 
$$
\omega=\sum_{(\ell_1,\nu_1),\cdots,(\ell_n,\nu_n)} a_{(\ell_1,\nu_1),\cdots,(\ell_n,\nu_n)}
(z_1,\dots,z_n)(d^{\ell_1}z_1)^{\nu_1}\cdots(d^{\ell_n}z_n)^{\nu_n}\leqno{(3)}
$$
on ${\bf C}^n$.

\medskip
Take a relatively compact open ball $B_0$ in ${\bf C}^n$ and a positive
number $b$ such that for every ball $W_k$ of radius $b$ in ${\bf C}^n$
there exists an element $\ell\in\Lambda$ such that the translate of
$W_k$ by $\ell$ is contained in $B_0$.
Since $\omega$ vanishes to order at least $q$ along $J_k(D)\cap(A\times W_k)$ 
it follows that there exists holomorphic jet differentials $\sigma_{j_0,\cdots,j_k}$ on
some open neighborhood $U$ of the topological closure of $B_0$ such that
$$
\omega=\sum_{j_0+\cdots+j_k=q}\sigma_{j_0,\cdots,j_k}\prod_{\nu=0}^k(d^\nu\theta_D)^{j_\nu}
$$ 
on $U\times W_k$.  Here we use $\omega$ to denote also its pullback to the universal
covering ${\bf C}^n$ of $A$. 
Clearly locally we can choose $\sigma_{j_0,\cdots,j_k}$
so that the order and the weight of each term
$\sigma_{j_0,\cdots,j_k}\prod_{\nu=0}^k(d^\nu\theta_D)^{j_\nu}$ are
respectively equal to those of $\omega$.
We rewrite the above identity on $U\times W_k$ in the form
$$
\omega=\theta_D^q\sum_{j_0+\cdots+j_k=q}\sigma_{j_0,\cdots,j_k}
\prod_{\nu=1}^k\left({d^\nu\theta_D\over{\theta_D}}\right)^{j_\nu}\leqno{(4)}
$$
on $U\times W_k$.  
In order to extend the identity (4) from $U\times W_k$ to identities on the translates
of $U\times W_k$ by elements of $\Lambda$, we rewrite (4) as
$$
\exp(-{\delta\over 2}\varphi)\omega=\exp(-{q\over 2}\varphi)
\theta_D^q\sum_{j_0+\cdots+j_k=q}{\hat\sigma}_{j_0,\cdots,j_k}
\prod_{\nu=1}^k\left({d^\nu\theta_D\over{\theta_D}}\right)^{j_\nu}\leqno{(5)}
$$
on $U$, where
$$
{\hat\sigma}_{j_0,\cdots,j_k}=\exp\left({{q-\delta}\over 2}\varphi\right)\sigma_{j_0,\cdots,j_k}.
$$
For $z\in{\bf C}^n$ and $\ell\in\Lambda$ there is an affine transformation
of ${\bf C}^n$ to itself given by
$z\mapsto A_\ell(z)+B_\ell$ such that
$$
\theta_D(z+\ell)=\exp(A_\ell(z)+B_\ell)\theta_D(z).
$$ 
From it we obtain
$$
d^p\theta_D(z+\ell)=\sum_{j=0}^p{p\choose j}
d^{p-j}\exp(A_\ell(z)+B_\ell)d^j\theta_D(z)
$$
and
$$
{d^p\theta_D(z+\ell)\over{\theta_D(z+\ell)}}=\sum_{j=0}^p{p\choose j}
\left({d^{p-j}\exp(A_\ell(z)+B_\ell)\over{\exp(A_\ell(z)+B_\ell)}}\right)
\left({d^j\theta_D(z)\over{\theta_D(z)}}\right).\leqno{(6)}
$$

\medskip
Both $\left|\exp(-{\delta\over 2}\varphi)\omega\right|$
and $\left|\exp(-{q\over 2}\varphi)\theta_D^q\right|$ are invariant under the action
of $\Lambda$.  Hence for $\ell\in\Lambda$ there exists a complex number
$c_{\ell}$ of absolute value $1$ such that
$$
c_\ell\left({\exp(-{\delta\over 2}\varphi)\omega
\over{\exp(-{q\over 2}\varphi)\theta_D^q}}\right)(z+\ell)
=\left({\exp(-{\delta\over 2}\varphi)\omega
\over{\exp(-{q\over 2}\varphi)\theta_D^q}}\right)(z).
$$

For any $z\in{\bf C}^n$ there exists $\ell=\ell(z)\in\Lambda$
such that the ball $B$ of radius $b$ centered at $z-\ell$ is contained in
$B_0$.
From
$$
\displaylines{
\left(\exp(-{\delta\over 2}\varphi)\omega\right)(z-\ell)=
\hfill\cr\hfill
\left(\exp(-{q\over 2}\varphi)
\theta_D^q\right)(z-\ell)
\sum_{j_0+\cdots+j_k=q}{\hat\sigma}_{j_0,\cdots,j_k}(z-\ell)
\prod_{\nu=1}^k\left({d^\nu\theta_D\over{\theta_D}}\right)^{j_\nu}(z-\ell)
\cr}
$$
on $B\times W_k$ and (6) we obtain
$$
\displaylines{
c_\ell\left(\exp(-{\delta\over 2}\varphi)\omega\right)(z)=
\hfill\cr\hfill
\left(\exp(-{q\over 2}\varphi)
\theta_D^q\right)(z)
\sum_{j_0+\cdots+j_k=q}P_{j_0,\cdots,j_k}
\prod_{\nu=1}^k\left({d^\nu\theta_D\over{\theta_D}}\right)^{j_\nu}(z)
\cr}
$$
on $(B+\ell)\times W_k$, where 
$$
P_{j_0,\cdots,j_k}^{(\lambda)}=P_{j_0,\cdots,j_k}^{(\lambda)}\left(
\{{\hat\sigma}_{p_0,\cdots,p_k}^{(\lambda)}(z-\ell)\}_{p_0+\cdots+p_k=q},
\{d^\nu(A_{-\ell}z+B_{-\ell})\}_{0\leq\nu\leq k}\right)
$$
is a polynomial in the variables 
$$
\{{\hat\sigma}_{p_0,\cdots,p_k}^{(\lambda)}\}_{p_0+\cdots+p_k=q},
\{d^\nu(A_{-\ell}z+B_{-\ell})\}_{0\leq\nu\leq k}
$$
whose coefficients are complex numbers.
Since $A_\ell$ is linear in $\ell$ and $B_\ell$ is a polynomial of degree
at most $1$ in $\ell$, 
it follows that
$$
\left|f^*d^{\nu}(A_{-l}z+B_{-l})\right|
\leq\left(1+\sum_{0\le\nu\le k,1\le j\le n}|f^*d^{\nu}z_j|^2\right).
$$
Hence
$$
f^*\left|\exp\left(-{\delta\over 2}\varphi\right)\omega\right|\leq\leqno{(7)}
$$
$$
C f^*\left|\exp\left(-{q\over 2}\varphi\right)\theta_D^q\right|
\left(1+\sum_{0\leq\nu\leq k,1\leq j\leq n}|f^*d^\nu z_j|\right)^N
\sum_{j=0}^k\left|f^*\left({d^j\theta_D\over{\theta_D}}\right)\right|^q
$$
on ${\bf C}$ where $N$ is a positive integer depending on $k$ and $q$ and
$C$ is a positive constant.

\medskip
Let ${\cal A}_r$ denote the operator which, when applied to a function, averages the
function over the circle of radius $r$ in ${\bf C}$ centered at the origin.
For a meromorphic function $g$ on $\bf C$ we denote by $T(r,g)$ 
the characteristic function of $F$ which is given by
$$
T(r,F)={\cal A}_r(\log^+|F|)+\int^r_{\rho =0}n(r,F,\infty){d\rho\over\rho},
$$
where $n(r,F,\infty)$ is the number of poles of $F$ with multiplicity
counted in the open disk of radius $r$ in $\bf C$ centered at the origin.

\medskip
To compute the defect for the map $f$ and the divisor $D$ we have to consider
$$
{\cal A}_r\left(\log^+{1\over{f^*\left|\exp(-{q\over 2}\varphi)\theta_D^q\right|}}\right)\leqno{(8)}
$$
which by (7) is dominated by
$$
{\cal A}_r\left(
\log^+{1\over{f^*\left|\exp(-{\delta\over 2}\varphi)\omega\right|}}
\right)
+\leqno{(9)}
$$
$$
{\cal A}_r\left(\log^+\left(
\left(1+\sum_{0\leq\nu\leq k,1\leq j\leq n}|f^*d^\nu z_j|\right)^N
\sum_{j=0}^k\left|f^*\left({d^j\theta_D\over{\theta_D}}\right)\right|^q\right)\right)
+O(1).
$$
Here $O(1)$ means the standard Landau symbol for order comparison.
We handle the first term of (9) as follows.
$$
{\cal A}_r\left(
\log^+{1\over{f^*\left|\exp(-{\delta\over 2}\varphi)\omega\right|}}
\right)
\leq{\cal A}_r\left({\delta\over 2}\varphi\circ f\right)
+{\cal A}_r\left(\log^+{1\over{\left|f^*\omega\right|}}\right)\leqno{(10)}
$$
$$
\leq{\cal A}_r\left({\delta\over 2}\varphi\circ f\right)
+T\left(r,{1\over{f^*\omega}}\right)
$$
$$
\leq{\cal A}_r\left({\delta\over 2}\varphi\circ f\right)
+T(r,f^*\omega)+O(1).
$$
Here according to our convention $f^*\omega$ is regarded as a function on
$\bf C$.  For (10) we have used the First Main Theorem of Nevanlinna
that 
$$
T\left(r,{1\over{f^*\omega}}\right)=T(r,f^*\omega)+O(1).
$$
The positive $(1,1)$-form 
${\sqrt{-1}\over{2\pi}}\partial{\overline\partial}\varphi$ is the curvature
form for the line bundle $L_D$ with the Hermitian metric $\exp(-\varphi)$.
We denote by $T(r,f,L_D)$ the characteristic function of $f$ with respect to
the Hermitian line bundle $L_D$ which is given by
$$
T(r,f,L_D)=\int^r_{\rho =0}{d\rho\over\rho}
\int_{\{|\zeta|<\rho\}}f^*{\sqrt{-1}\over{2\pi}}\partial{\overline\partial}\varphi
$$
which by Green's theorem equals
${\cal A}_r\left({\delta\over 2}\varphi\circ f\right)+O(1)$.

From (2) and (3) we conclude that
$$
T(r,f^*\omega)\leq\delta T(r,f,L_D)+O\left(\log T(r,f,L_D)\right).\leqno{(11)}
$$
Here we have used the fact that
$$
T(r,f^*d^\nu z_j)=O\left(\log T(r,f,L_D)\right),\leqno{(12)}
$$ 
which is a consequence of the logarithmic derivative lemma,
when $f^*d^\nu z_j$ is regarded as a meromorphic function on $\bf C$.
The second term of (9) satisfies
$$
{\cal A}_r\left(\log^+\left(
\left(1+\sum_{0\leq\nu\leq k,1\leq j\leq n}|f^*d^\nu z_j|\right)^N
\sum_{j=0}^k\left|f^*\left({d^j\theta_D\over{\theta_D}}\right)\right|^q\right)\right)\leqno{(13)}
$$
$$
=
O\left(\log T(r,f,L_D)\right)
$$ 
because of (12) and because of
$$
T\left(r,f^*\left({d^j\theta_D\over{\theta_D}}\right)\right)=O\left(\log T(r,f,L_D)\right),
$$
which is a consequence of the logarithmic derivative lemma,
when $f^*\left({d^j\theta_D\over{\theta_D}}\right)$ is regarded as a meromorphic function on $\bf C$.
Finally from the domination of (8) by (9) and from (10) and (11) and (13)
we conclude that
$$
{\cal A}_r\left(\log^+{1\over{f^*\left|\exp(-{q\over 2}\varphi)\theta_D^q\right|}}\right)
=2\delta T(r,f,L_D)+O\left(\log T(r,f,L_D)\right).\leqno{(14)}
$$
We denote by $m(r,f,D)$ the proximity function 
for the map $f$ and the divisor $D$ which is defined, up to a bounded term,
by
$$
m(r,f,D)={\cal A}_r\left(\log{1\over{|s_D|}}\right),
$$
where $s_D$ is the canonical section of $L_D$ whose divisor is $D$ and
$|s_D|$ is the pointwise norm of $s_D$ with respect to a Hermitian metric
of $L_D$.  From this definition of $m(r,f,D)$ we have
$$
{\cal A}_r\left(\log^+{1\over{f^*\left|\exp(-{q\over 2}\varphi)\theta_D^q\right|}}\right)
=q\ m(r,f,D) + O(1).
$$
For any given $\epsilon>0$ we can choose $\delta$ and $q$ so that ${\delta\over q}<\epsilon$.
We denote by $\delta(f,D)$ the defect for the map $f$ and the divisor $D$
which is defined by
$$
\delta(f,D)={\hbox{lim inf}}_{r\rightarrow\infty}{m(r,f,D)\over{T(r,f,L_D)}}.
$$
It follows from (14) that $\delta(f,D)<{2\delta\over q}<2\epsilon$.  
From the arbitrariness of the positive number $\epsilon$
we conclude that the defect for an ample divisor in an abelian variety
is $0$.

\bigskip
\noindent
{\bf \S2. The General Schwarz Lemma for the Holomorphic Case}.  

\medskip
We introduce the following notation.
For a function or a (1,1)-form $\eta$ let 
$$
{\cal I}_r(\eta)= 
\int^r_{\rho =0}{d\rho\over\rho}\int_{\zeta\in{\bf C},|\zeta|<\rho}\eta,
$$  
where $\zeta$ is the coordinate of $\bf C$.  Green's theorem gives
$$
{\cal I}_r\left({1\over\pi}\partial_\zeta\partial 
_{\overline\zeta}
g\right)={1\over 2}{\cal A}_r(g)-{1\over 2} g(0)
$$
for a function $g$.  For a complex manifold $M$ and a positive definite
smooth $(1,1)$-form $\theta$ on $M$ and for a holomorphic map $f:{\bf C}
\rightarrow M$, we define the characteristic function of $f$ with respect
to $\theta$ as
$$
T(r,f,\theta)={\cal I}_r(f^*\theta).
$$
For another positive definite smooth $(1,1)$-form $\theta^\prime$,
$T(r,f,\theta)\leq CT(r,f,\theta^\prime)$ for some constant $C$ depending
only on $\theta$ and $\theta^\prime$ and independent of $f$.  
When it does not matter which positive definite smooth $(1,1)$-form 
$\theta$ is used, we also denote $T(r,f,\theta)$ simply by $T(r,f)$.
If $M$ is the complex projective line ${\bf P}_1$ and $f:{\bf C}
\rightarrow{\bf P}_1$ is represented by a meromorphic function $F$
on $\bf C$, then
$$
T(r,f,\theta)=T(r,F)+O(1)
$$
when $\theta$ is the Fubini-Study form on ${\bf P}_1$.

\medskip
For a holomorphic line bundle $L$ with a smooth 
Hermitian metric $e^{-\varphi}$
over a compact complex manifold $M$ and an $L$-valued meromorphic
$k$-jet differential $\omega$, we define the pointwise norm
$$
|\omega|_L=\left(e^{-\varphi}\omega\overline\omega\right)^{1\over 2}.
$$
Locally the pointwise norm is the absolute value of a meromorphic jet 
differential.  The definition of $|\omega|_L$ does not involve any
metric of the tangent bundle of $M$.
At a point $P$ of $M$ the pointwise norm
$|\omega|_L$ is not a scalar.  However, if a $k$-jet $\xi$ is given at $P$
not in the pole set of $\omega$
(for example, the $k$-jet defined by a holomorphic map from an open subset
of ${\bf C}$ whose image contains $P$),
the value of $|\omega|_L$ at $\xi$ is a nonnegative number.  The pointwise
norm $|\omega|_L$ depends not just on the line bundle $L$ but also on
the metric $e^{-\varphi}$ of $L$, 
but we will simply use the notation $|\omega|_L$
and suppress the metric $e^{-\varphi}$ if there is no confusion.  If $L$
is the line bundle associated to a divisor $D$, we also use $D$ to denote
the line bundle $L$ and denote $|\omega|_L$ by $|\omega|_D$ if there is no 
confusion.  We use the additive notation for tensor products of line
bundles when the expression involves using a divisor to denote the line
bundle associated to it.  If $L$ is the trivial line bundle, we denote
$|\omega|_L$ simply by $|\omega|$.

\medskip
For a meromorphic jet differential $\eta$ on $M$ there is an ample divisor 
$D$ with canonical section $s_D$ such that $s_D^k\eta$ is a $kD$-valued
holomorphic jet differential.  The pointwise norm $|\eta|$ for the 
meromorphic jet differential $\eta$ is not the same as the pointwise norm
$|s_D^k\eta|_{kD}$ of the $kD$-valued holomorphic jet differential
$s_D^k\eta$ when $kD$ is given a smooth Hermitian metric
(though a meromorphic
jet differential $\eta$ can equivalently be regarded a $kD$-valued 
holomorphic jet differential $s_D^k\eta$).  The value
of $|\eta|$ can blow up when evaluated at a smooth field of jets, but
the value of $|s_D^k\eta|_{kD}$ at a smooth field of jets is smooth.  The 
distinction between the
pointwise norm of a meromorphic jet differential and the pointwise norm
of the line-\-bundle-\-valued
holomorphic jet differential associated to it is
crucial in our proof of the general Schwarz lemma.

\medskip
\noindent
{\it Lemma 3}.  Let $M$ be a compact complex manifold, $E$ be a holomorphic
line bundle over $M$ with hermitian metric $e^{-\psi}$ along the fibers of
$E$, and curvature $\theta_\psi=\partial\overline\partial\psi$.
Let ${\cal D}_\psi$ denote covariant differentiation with respect to the
connection from the metric $e^{-\psi}$ of $E$.  
Let $\omega$ be an $E$-valued holomorphic jet differential over $M$.
Then
$$
\partial\overline\partial\log\left(1+|\omega|_E^2\right)
={{|{\cal D}_\psi\omega|_E^2}\over{(1+|\omega|_E^2)^2}}
-{{\theta_\psi|\omega|_E^2}\over{1+|\omega|_E^2}},\leqno{(15)}
$$
where all differentials are taken in the category
of jet differentials instead of in the category of differential forms
({\it i.e.}, the differentials are symmetric instead of alternating). 
Moreover, if $f:{\bf C}\rightarrow M$ is a holomorphic map, then
$$
{\cal I}_r\left(
f^*{{|{\cal D}_\psi\omega|_E^2}\over{(1+|\omega|_E^2)^2}}
\right)
\leq
{\cal A}_r\left(
f^*\log\left(1+|\omega|_E^2\right)
\right)+O(T(r,f,\theta_\psi)+\log r).
$$

\medskip
\noindent
{\it Proof}.  The identity (15)
follows from straightforward differentiation.  Note that,
if the order of the jet differential $\omega$ is $k$, 
then ${\cal D}_\psi\omega$ is a jet differential of order $k+1$ which in general
is not holomorphic.  The last inequality follows from (15),
Green's theorem, and 
$$
{{|\omega|_E^2}\over{1+|\omega|_E^2}}\leq 1.
$$
Q.E.D.

\medskip
\noindent
{\it Definition}.  Let $M$ be a compact complex manifold and $D$ be an
ample divisor with canonical section $s_D$.  By a {\it meromorphic
jet differential constructed from functions with poles along $D$} we  
mean a meromorphic jet differential $\eta$ on $M$ which is a polynomial of
the variables $d^\nu\left({s_\lambda\over s_D}\right)$ with constant
coefficients, where $\nu$ is a nonnegative integer and $s_\lambda$ is
a holomorphic section of the line bundle associated to $D$.  By a 
{\it pole-factor} of $\eta$ we mean $s_D^k$ so that $s_D^k\eta$ is a 
holomorphic jet differential.

\medskip
\noindent
{\it Lemma 4.} Let $k$ and $m$ be positive integers.  
Let $M$ be a compact complex projective algebraic manifold,
$E$ be a holomorphic line bundle over $M$ and $\omega$ be an $E$-valued
holomorphic $k$-jet differential of weight $m$ on $M$.  Then there
exist 

\medskip
\noindent
(i) an ample line bundle $F$ with holomorphic sections $s_{D_\nu}$ whose
divisor is a nonsingular ample divisor $D_\nu$ ($1\leq\nu\leq N$),

\medskip
\noindent
(ii) meromorphic $k$-jet differentials $\eta_\nu$ weight $m$
constructed by functions with
poles along $D_\nu$ (which automatically
admits $s_{D_\nu}^{2m}$ as a pole-factor) for $1\leq\nu\leq N$, and

\medskip
\noindent
(iii) an ample line bundle $L$ over $M$ and holomorphic sections
$t_{j,\nu}$ of $L+E$ over $M$ and holomorphic sections
$t_j$ of $L+2mF$ ($1\leq j\leq J, 1\leq\nu\leq N$)

\medskip
\noindent
such that

\medskip
\noindent
($\alpha$) $t_1,\cdots,t_J$ have no common zeroes in $M$, and

\medskip
\noindent
($\beta$) $t_j\omega=\sum_{\nu=1}^Nt_{j,\nu}\left(s_{D_\nu}^{2m}\eta_\nu\right)$
on $M$ for $1\leq j\leq J$.

\medskip
In particular,
$$
|\omega|_E\leq 
C\sum_{j=1}^J\sum_{\nu=1}^N\left|s_{D_\nu}^{2m}\eta_\nu\right|_{2mF}
$$
on $M$ for some constant $C$.

\medskip
\noindent
{\it Proof}.  It follows from the standard application to a short exact 
sequence of the vanishing theorem for the positive 
dimensional cohomology over a compact complex manifold with coefficients 
in the tensor product of a holomorphic vector bundle with a sufficiently 
high power of an ample line bundle.  Q.E.D.

\medskip
\noindent
{\it Proposition 1.} Let $M$ be a compact projective-algebraic complex
manifold of complex dimension $n$ and $f:{\bf C}\rightarrow M$ be
a nonconstant holomorphic map.  Let $E$ be a holomorphic line bundle with
Hermitian metric $e^{-\psi}$ along its fibers.  Let $\omega$ be
an $E$-valued holomorphic $k$-jet differential on $M$ of weight $m$.  
Then there exists a positive number $\epsilon_{k,m}$ depending only on $k$
and $m$ such that for $0<\epsilon\leq\epsilon_{k,m}$ one has
$$
{\cal I}_r\left(\left(f^*|\omega|^2_E\right)^\epsilon\right)
=O\left(T(r,f)+\log r\right)\quad \|.
$$

\medskip
\noindent
{\it Remark}.
Note that in this paper the pullback $f^*\omega$ means the value
of $\omega$ at $d^kf$ with respect to a local trivialization of $E$ which
means that we take only the part of $f^*\omega$ which is equal to a scalar
function times $(d\zeta)^m$.  As customary in value distribution theory,
we use the notation $\|$ at the end of an equation or an inequality to 
mean that the statement holds outside an open set whose harmonic measure is
finite.  Proposition 1 is the step in the proof of Theorem 1 which 
corresponds to using Cauchy's integral formula for derivatives to dominate
higher order derivatives at a point by first-order derivative on a circle
centered at that point.  The condition restricting the validity of the
statement to outside an open set with finite harmonic measure is used to
take care of the problem posed by the need for a circle in the use of
Cauchy's integral formula for derivatives.

\medskip
\noindent
{\it Proof of Proposition 1}.  By Lemma 4 it suffices to prove the special case that
$E$ is equal to $2mD$ for some ample divisor $D$ of $M$ with canonical 
section $s_D$ and $\omega=s_D^{2m}\eta$ for some meromorphic $k$-jet
differential $\eta$ of weight $m$ constructed from functions with
poles along $D$.  From the definition of meromorphic jet
differential constructed from functions with
poles along $D$ we can further assume that
$$
\eta=\prod_{1\leq\nu\leq k,1\leq\lambda\leq n}
\left(d^\nu\left({s_\lambda\over s_D}\right)\right)^{\ell_{\nu,\lambda}},
$$
where $s_\lambda$ ($1\leq\lambda\leq n$) are holomorphic sections of
$D$ over $M$ and 
$$
\sum_{1\leq\nu\leq k,1\leq\lambda\leq n}\nu\ell_{\nu,\lambda}=m.
$$
By using the trivial inequality
$$
\displaylines{
\left|\prod_{1\leq\nu\leq k,1\leq\lambda\leq n}
\left(d^\nu\left({s_\lambda\over s_D}\right)\right)^{\ell_{\nu,\lambda}}\right|
\leq
\left(
\max_{1\leq\nu\leq k,1\leq\lambda\leq n}
\left|d^\nu\left({s_\lambda\over s_D}\right)\right|\right)
^{\sum_{1\leq i\leq k,1\leq j\leq n}\ell_{i,j}}\cr
\leq
\sum_{\nu=1}^k\sum_{\lambda=1}^n
\left|d^\nu\left({s_\lambda\over s_D}\right)\right|
^{\sum_{1\leq i\leq k,1\leq j\leq n}\ell_{i,j}},\cr
}
$$
we need only prove the special case $\omega=s_D^{\nu+1}d^\nu\left({s_\lambda\over s_D}\right)$.
We prove the special case by induction on $\nu$.
The case of $\nu=1$ is
clear.  Assume $\nu>1$. Let 
$$
\eta^\prime=d^{\nu-1}\left({s_\lambda\over s_D}\right).
$$
Then 
$d\eta^\prime=\eta$.
We use the symbols $C_\nu$ ($1\leq\nu\leq 7$)
to denote constants.
Let $e^{-\psi}$ be a smooth Hermitian metric for $D$
and $\eta_0$ be a smooth positive definite $(1,1)$-form on $M$.
Then
$$
\left|s_D{\cal D}_{\nu\psi}\left(s_D^\nu\eta^\prime\right)
-s_D^{\nu+1}\eta\right|_{(\nu+1)D}
\leq C_1\left|s_D^\nu\eta^\prime\right|_{\nu D}(\eta_0)^{1\over 2},\leqno{(16)}
$$
where ${\cal D}_{\nu\psi}$ denote the covariant differentiation for sections
of the line bundle $\nu D$ with respect to the metric $e^{-\nu\psi}$.
By Lemma 3 and (16) we have
$$
{\cal I}_r\left(
f^*{{
\left|s_D^{\nu+1}\eta\right|_{(\nu+1)D}^2
}\over{(1+|s_D^\nu\eta^\prime|^2)^2}}
\right)
\leq C_2\cdot
{\cal A}_r\left(
f^*\log\left(1+|s_D^\nu\eta^\prime|_{\nu D}^2\right)
\right)\leqno{(17)}
$$
$$
+O(T(r,f)+\log r).
$$
Take $0<\epsilon\leq{1\over 2}\min\left(1,\epsilon_{\nu,\nu}\right)$.
From (17) it follows that
$$
{\cal I}_r\left(\left(
f^*
\left|s_D^{\nu+1}\eta\right|_{(\nu+1)D}^2
\right)^\epsilon\right)
=
{\cal I}_r\left(\left(
f^*{{
\left|s_D^{\nu+1}\eta\right|_{(\nu+1)D}^2
}\over{(1+|s_D^\nu\eta^\prime|_{\nu D}^2)^2}}
f^*(1+|s_D^\nu\eta^\prime|_{\nu D}^2)^2
\right)^\epsilon\right)
\leqno{(18)}
$$
$$
\displaylines{
\leq
C_3\left(
{\cal I}_r(1)+
{\cal I}_r\left(
f^*{{
\left|s_D^{\nu+1}\eta\right|_{(\nu+1)D}^2
}\over{(1+|s_D^\nu\eta^\prime|_{\nu D}^2)^2}}
\right)+
{\cal I}_r\left(\left(
f^*(1+|s_D^\nu\eta^\prime|_{\nu D}^2)^2\right)^{2\epsilon}\right)
\right)
\cr
\leq
C_4\left(\left(
{\cal A}_r\left(
f^*\log\left(1+|s_D^\nu\eta^\prime|_{\nu D}^2\right)
\right)
\right)
+{\cal I}_r\left(\left(
f^*(1+|s_D^\nu\eta^\prime|_{\nu D}^2)^2\right)^{2\epsilon}\right)\right)
\cr
+O(T(r,f)+\log r)
\cr
\leq
C_5{\cal A}_r\left(
f^*\log\left(1+|s_D^\nu\eta^\prime|_{\nu D}^2\right)
\right)+O(T(r,f)+\log r)\quad\|,
\cr}
$$
where for the last inequality the induction hypothesis is used.
The standard techniques in value distribution theory of using
the so-called Calculus Lemma and the concavity of 
the logarithmic function gives
$$
{\cal A}_r(\log g)\leq
C_6\left(\log r+\log{\cal I}_r(g)\right)
\quad \|\leqno{(19)}
$$
for any smooth positive function $g$ on $\bf C$.
Putting this into (18) yields
$$
\displaylines{
{\cal I}_r\left(\left(
f^*
\left|s_D^{\nu+1}\eta\right|_{(\nu+1)D}^2
\right)^\epsilon\right)\cr
\leq
C_7\log{\cal I}_r\left(
f^*\left(1+|s_D^\nu\eta^\prime|_{\nu D}^2\right)^{2\epsilon}
\right)+O(T(r,f)+\log r)\ \ \|,\cr
}
$$
when we use
$$
g=f^*\left(1+|s_D^\nu\eta^\prime|_{\nu D}^2\right)^{2\epsilon}.
$$
Thus
$$
{\cal I}_r\left(
f^*{{
\left|s_D^{\nu+1}\eta\right|_{(\nu+1)D}^2
}\over{(1+|s_D^\nu\eta^\prime|_{\nu D}^2)^2}}
\right)\leq O(\log T(r,f)+\log r)\quad\|
$$
by induction hypothesis and the induction argument is complete. Q.E.D.

\medskip
\noindent
{\it Proof of Theorem 1.} Let $s_D$ be the canonical section of the line bundle 
$L_D$ associated
to the divisor $D$.  Let $e^{-\psi}$ be a Hermitian metric of $L_D$ 
whose curvature is a positive definite smooth $(1,1)$-form $\theta_D$ on
$X$ (in the sense that $\theta_D$ is the restriction of some smooth
positive definite $(1,1)$-form on ${\bf P}_N$).  Suppose $f^*\omega$ is
not identically zero on $\bf C$ and we are going to derive a contradiction.
From
$$
\partial\overline\partial\log
\left(
e^\psi
\left(
{\omega\over{s_D}}
\right) 
\overline{\left(
{\omega\over{s_D}}
\right)}
\right)
\geq
\theta_D
$$
we conclude from Green's theorem that
$$
{\cal A}_r\left(
\log
f^*\left|{\omega\over{s_D}}\right|_{-D}^2
\right)
\geq
{\cal I}_r\left(f^*\theta_D\right)-O(1).\leqno{(20)}
$$
It follows from (19) that
$$
{\cal A}_r\left(\log 
f^*\left|{\omega\over{s_D}}\right|_{-D}^\epsilon
\right)\leq
C\left(\log r+\log{\cal I}_r\left(f^*
\left|{\omega\over{s_D}}\right|_{-D}^\epsilon
\right)\right)
\quad \|
$$
for any $\epsilon>0$ and for some constant $C$ depending on $\epsilon$.
Let $H_N$ be the hyperplane section line bundle of ${\bf P}_N$.  There
exist some positive integer $\ell$ and a number of global holomorphic sections
$\sigma_1,\cdots,\sigma_q$ of $H_N^{\otimes\ell}$ over ${\bf P}_N$ without
common zeroes so that
$\sigma_j\left({\omega\over s_D}\right)$ can be extended to a global
$H_N^{\otimes\ell}$-valued holomorphic $k$-jet differential over ${\bf P}_N$
for $1\leq j\leq q$.
By Proposition 1 applied to 
$\sigma_j\left({\omega\over s_D}\right)$ on ${\bf P}_N$ for $1\leq j\leq q$, 
we get
$$
{\cal I}_r\left(
f^*\left|{\omega\over{s_D}}\right|_{-D}^\epsilon
\right)=
O(T(r,f)+\log r)\quad\|\leqno{(21)}
$$
for some $\epsilon>0$
and consequently
$$
{\cal A}_r\left(\log 
f^*\left|{\omega\over{s_D}}\right|_{-D}^\epsilon
\right)=
O(\log T(r,f)+\log r)\quad\|.
$$
Combining with (20), we obtain
$$
{\cal I}_r\left(f^*\theta_D\right)
\leq
O(\log T(r,f)+\log r)\quad\|,
$$
which implies that $T(r,f)$ is of the order $\log r$.  From (21) it
follows that
$$
{\cal I}_r\left(f^*|\omega|^\epsilon\right)=O(\log r)\quad\|.\leqno{(22)}
$$
On the other hand, from the subharmonicity of $f^*|\omega|^\epsilon$
we conclude that the growth of ${\cal I}_r\left(f^*|\omega|^\epsilon\right)$
is at least $r^2$, contradicting (22). Q.E.D.

\bigskip
\noindent
{\bf \S3. The General Schwarz Lemma for the Log-Pole Case}.  

\medskip
Theorem 2 is proved by a modification of the arguments in the
proof of Theorem 1.  We only present here the necessary modifications.
For the modification of Proposition 1, we assume that $X$ is nonsingular and
let $M=X$.  We let $t_j$ be the canonical section of the line bundle
over $M$ associated to the divisor $Z_j$ ($1\leq j\leq q$).  We assume that
the image of the holomorphic map $f$ is disjoint from $\cup_{j=1}^qZ_j$.
Fix a
smooth metric for the line bundle $Z_j$ so that the norm of $|t_j|_{Z_j}<1$
on $M$.  For any positive number $A>e$ let 
$\tau_{j,A}=\log\left({A\over|t_j|_{Z_j}^2}\right)$.  The modified Proposition 1
states that, for any $E$-valued meromorphic $k$-jet differential $\omega$ of
weight $m$ on $M$ with at most log-pole singularity along $\cup_{j=1}^qZ_j$,
there exist positive numbers $\epsilon_{k,m}$ and $A_{k,m}$ and a positive
integer $a_{k,m}$ such that
$$
{\cal I}_r\left(\left(f^*{|\omega|_E\over
\prod_{j=1}^q\tau_{j,A}^a}\right)^\epsilon\right)
=O\left(T(r,f)+\log r\right)\quad \|
$$
for $0<\epsilon\leq \epsilon_{k,m}$, $A>A_{k,m}$, and $a\geq a_{k,m}$.

\medskip
\noindent
{\it Remark.} Heuristically speaking, the factor 
$\prod_{j=1}^q\tau_{j,A}^a$ in the statement of the modified
Proposition 1 is due to the restriction placed
by the log-pole $\cup_{j=1}^q Z_q$ on the radius of the circle $\Gamma$
used in 
the domination of the higher
order derivatives at the center of $\Gamma$
by the first-order derivative on $\Gamma$
by Cauchy's integral formula for derivatives.

\medskip
Lemma 3 holds when for $\omega$ with log-pole singularities as long as the
image of the map $f$ is disjoint from the log-pole.  We use
$C_j$ ($1\leq j\leq 6$) to denote constants.  From the last
inequality in Lemma 3 it follows that
$$
\displaylines{
{\cal I}_r\left(\left(
f^*{{|{\cal D}_\psi\omega|_E^2}\over\prod_{j=1}^q\tau_{j,A}^a}
\right)^\epsilon\right)\leq
C_1\left\{
{\cal I}_r\left(
f^*{{|{\cal D}_\psi\omega|_E^2}\over{(1+|\omega|_E^2)^2}}
\right)+
{\cal I}_r\left(\left(
f^*{{1+|\omega|_E^2}\over\prod_{j=1}^q\tau_{j,A}^a}
\right)^{2\epsilon}\right)
\right\}
\cr
\leq
C_2\left\{
{\cal A}_r\left(
f^*\log\left(1+|\omega|_E^2\right)
\right)+
{\cal I}_r\left(\left(
f^*{|\omega|_E^2\over\prod_{j=1}^q\tau_{j,A}^a}
\right)^{2\epsilon}\right)
\right\}\cr
+O(T(r,f,\theta_\psi)+\log r)\cr
\leq
C_3\Bigg\{
{\cal A}_r\left(
f^*\log\left(1+{|\omega|_E^2\over\prod_{j=1}^q\tau_{j,A}^a}\right)
\right)+
{\cal A}_r\left(
f^*\log\prod_{j=1}^q\tau_{j,A}^a
\right)+\cr
{\cal I}_r\left(\left(
f^*{|\omega|_E^2\over\prod_{j=1}^q\tau_{j,A}^a}
\right)^{2\epsilon}\right)
\Bigg\}
+O(T(r,f,\theta_\psi)+\log r)\cr
\leq
C_4\left\{
\sum_{j=1}^q\log{\cal A}_r\left(
f^*\tau_{j,A}
\right)+
{\cal I}_r\left(\left(
f^*{|\omega|_E^2\over\prod_{j=1}^q\tau_{j,A}^a}
\right)^{2\epsilon}\right)
\right\}\cr
+O(T(r,f,\theta_\psi)+\log r)\quad\|.\cr
}
$$
Using the domination of the proximity function 
${\cal A}_r\left(f^*\tau_{j,A}\right)$
by the characteristic function $T(r,f,\theta_{Z_j})+O(1)$ 
(where $\theta_{Z_j}$
is the curvature of the metric of the line bundle $Z_j$), we conclude that
$$
\displaylines{
{\cal I}_r\left(\left(
f^*{{|{\cal D}_\psi\omega|_E^2}\over\prod_{j=1}^q\tau_{j,A}^a}
\right)^\epsilon\right)\cr
\leq
C_5\left\{
{\cal I}_r\left(\left(
f^*{|\omega|_E^2\over\prod_{j=1}^q\tau_{j,A}^a}
\right)^{2\epsilon}\right)
\right\}
+O(T(r,f,\theta_\psi)+\log r)\quad\|.\cr}
$$

\medskip
For the proof of the modified Proposition 1, it suffices to consider the
case where 
$\omega$ is of the form
$vt^{\nu+1}d^\nu\log\left({s\over t}\right)$,
where 

\medskip
\noindent
(i) $\nu$ is a positive integer,

\medskip
\noindent
(ii) $s,t$ are global holomorphic sections of some ample line 
bundle $L$,

\medskip
\noindent
(iii) no $Z_j$ ($1\leq j\leq q$) is a branch of the zero-set of $t$,

\medskip
\noindent
(iv) $v$ is a global holomorphic section of some ample line 
bundle $L^\prime$,

\medskip
\noindent
(v) no $Z_j$ ($1\leq j\leq q$) is a branch of the zero-set of $v$, and

\medskip
\noindent
(vi) $\omega$ is holomorphic outside $\cup_{j=1}^qZ_j$.

\medskip
As in the proof of Proposition 1, one further reduces the proof of the 
modified Proposition 1 to the case where $\nu=1$.  We choose a smooth
metric $e^{-\psi}$ of $L$ so that its curvature $\theta_\psi$ is a 
positive definite $(1,1)$-form on $M$.  We replace $e^{-\psi}$ by
$Ae^{-\psi}$ for a sufficiently large positive number $A$ so that
$|s|_L<1$ and $-\log|s|_L^2>{1\over\delta}$ on $M$ for some positive
number $\delta<1$.  Then
$$
{\sqrt{-1}\over 2\pi}\partial\overline\partial\log
{\left|t^2d\left({s\over t}\right)\right|^2_{L^{\otimes 2}}
\over
|s|_L^2(\log|s|_L^2)^2
}
\geq
-(1+\delta)\theta_\psi-\hbox{div}\, s+
{
|{\cal D}_\psi s|^2
\over
|s|^2\log|s|^2
},
$$
where $\hbox{div}\, s$ is the divisor of $s$ regarded as a $(1,1)$-current.
Now use
$$
{
\left|t{\cal D}_\psi s-t^2d\left({s\over t}\right)\right|^2_{L^{\otimes 2}}
\over
|s|^2_L
}
=O(1)
$$
and the standard techniques in value distribution theory of Green's Theorem,
the Calculus Lemma, and the concavity of the logarithmic function to 
conclude that
$$
{\cal I}_r\left(
f^*{\left|t^2d\left({s\over t}\right)\right|^2_{L^{\otimes 2}}
\over
|s|_L^2(\log|s|_L^2)^2
}\right)=O(T(r,f)+\log r)\quad\|.
$$
Because of condition (vi), the section $v$ used as a factor in
$\omega$ takes care of all the branches
of the divisor of $s$ except those contained in the log-pole set
$\cup_{j=1}^q Z_j$.
This concludes the proof of the modified Proposition 1.

\medskip
Now we come to the proof of Theorem 2.  Assume that $f^*\omega$ is
not identically zero and we are going to derive a contradiction.
The same argument used in 
the proof of Theorem 1 gives us (20) from which we conclude
$$
\displaylines{
{\cal I}_r\left(f^*\theta_D\right)
\leq
C_1\left(
{\cal A}_r\left(
\log
f^*{\left|{\omega\over{s_D}}\right|_{-D}^2\over\prod_{j=1}^q\tau_{j,A}^a}
\right)
+
{\cal A}_r\left(
\log
f^*\prod_{j=1}^q\tau_{j,A}^a
\right)
\right)
+O(1)\cr
\leq
C_2\left(
\log{\cal I}_r\left(
f^*\left({\left|{\omega\over{s_D}}\right|_{-D}^2\over\prod_{j=1}^q\tau_{j,A}^a}
\right)^\epsilon\right)
+
\sum_{j=1}^q\log{\cal A}_r\left(
f^*\tau_{j,A}
\right)
\right)
+O(\log r)\cr
=
O(\log T(r,f)+O(\log r))\quad\|,\cr}
$$
where $0<\epsilon\leq\epsilon_{k,m}$.
This implies that $T(r,f)$ is of the order $\log r$.
From
$$
\partial\overline\partial\log
{1\over\tau_{j,A}}
\geq
{-2\over\log\left|{A\over t_j}\right|^2}
\theta_{Z_j}
$$
it follows that
$$
\log f^*\left({|\omega|^2\over \prod_{j=1}^q\tau_{j,A}^a}\right)
$$
is subharmonic when $A$ is greater than some constant
${\tilde A}(a)$ depending on $a$ and consequently
$$
f^*\left({|\omega|^2\over \prod_{j=1}^q\tau_{j,A}^a}\right)^\epsilon
$$
is subharmonic, from which we conclude that
the order of growth of 
$$
{\cal I}_r\left(
f^*\left({|\omega|^2\over \prod_{j=1}^q\tau_{j,A}^a}\right)^\epsilon
\right)
$$
is at least $r^2$ which is a contradiction.

\bigskip
\noindent
{\bf \S4. Construction of Hyperbolic Hypersurfaces}

\medskip
We now construct hyperbolic hypersurfaces of degree
$16(n-1)^2$ in ${\bf P}_n$ by using Borel's lemma and a simple
dimension counting argument for certain subvarieties of Grassmannians.
We introduce a generalized Borel lemma (Proposition 2)
which will be used in \S5 and
which implies as an corollary the usual Borel lemma (Proposition 3) used in
here in \S4.

\medskip
\noindent
{\it Proposition 2 (Generalized Borel Lemma).}
Let $g_j(x_0,\cdots,x_n)$ be a homogeneous polynomial of degree $\delta_j$
for $0\leq j\leq n$.  Suppose there exists a holomorphic map 
$f:{\bf C}\rightarrow{\bf P}_n$ so that its image lies in
$$
\sum_{j=0}^nx_j^{p-\delta_j}g_j(x_0,\cdots,x_n)=0
$$
and $p>(n+1)(n-1)+\sum_{j=0}^n\delta_j$.
Then there is a nontrivial linear relation among
$x_1^{p-\delta_1}g_1(x_0,\cdots,x_n),\cdots,x_n^{p-\delta_n}g_n(x_0,\cdots,x_n)$
on the image of $f$.

\medskip
\noindent
{\it Proof}.  Use the affine coordinates $z_j={x_j\over x_0}$ for
$1\leq j\leq n$.  Let 
$$
{\tilde g}_j(z_1,\cdots,z_n)=x_0^{-\delta_j}g_j(x_0,\cdots,x_n).
$$
From
$$
\sum_{j=0}^nx_j^{p-\delta_j}g_j(x_0,\cdots,x_n)=0
$$
we have the following relation between the two Wronskians
$$
W\left({\tilde g}_0,z_1^{p-\delta_1}{\tilde g}_1,\cdots,z_{n-1}^{p-\delta_{n-1}}{\tilde g}_{n-1}\right)
=(-1)^{n-1}
W\left({\tilde g}_0,z_2^{p-\delta_2}{\tilde g}_2,\cdots,z_n^{p-\delta_n}{\tilde g}_n\right).
$$
Rewrite the equation as
$$
\displaylines{
z_1^{p-\delta_1-n+1}
\left\{{1\over\prod_{j=1}^{n-1}z_j^{p-\delta_j-n+1}}
W\left({\tilde g}_0,z_1^{p-\delta_1}{\tilde g}_1,\cdots,z_{n-1}^{p-\delta_{n-1}}{\tilde g}_{n-1}\right)\right\}\cr
=(-1)^{n-1}
z_n^{p-\delta_n-n+1}
\left\{{1\over\prod_{j=2}^nz_j^{p-\delta_j-n+1}}
W\left({\tilde g}_0,z_2^{p-\delta_2}{\tilde g}_2,\cdots,z_n^{p-\delta_n}{\tilde g}_n\right)\right\}.\cr}
$$
Since
$$
{1\over\prod_{j=1}^{n-1}z_j^{p-\delta_j-n+1}}
W\left({\tilde g}_0,z_1^{p-\delta_1}{\tilde g}_1,\cdots,z_{n-1}^{p-\delta_{n-1}}{\tilde g}_{n-1}\right)
$$
and 
$$
{1\over 
\prod_{j=2}^nz_j^{p-\delta_j-n+1}
}
W\left({\tilde g}_0,z_2^{p-\delta_2}{\tilde g}_2,\cdots,z_n^{p-\delta_n}{\tilde g}_n\right).
$$
are both holomorphic on the affine part of the hypersurface, we conclude
that
$$
{1\over
\prod_{j=1}^{n-1}z_j^{p-\delta_j-n+1}
}
W\left({\tilde g}_0,z_1^{p-\delta_1}{\tilde g}_1,\cdots,z_{n-1}^{p-\delta_{n-1}}{\tilde g}_{n-1}\right)
$$
is divisible by $z_n^{p-\delta_n-n+1}$.  Thus
$$
{1\over
\prod_{j=1}^nz_j^{p-\delta_j-n+1}}
W\left({\tilde g}_0,z_1^{p-\delta_1}{\tilde g}_1,\cdots,z_{n-1}^{p-\delta_{n-1}}{\tilde g}_{n-1}\right)
$$
is holomorphic on the affine part of the hypersurface.  Now we look at the
infinity part.  We introduce the coordinates
$w_j={x_j\over x_n}$ for $0\leq j\leq{n-1}$ so that
$z_j={w_j\over w_0}$ for $1\leq j\leq{n-1}$ and $z_n={1\over w_0}$.
We have
$$
\displaylines{
W\left({\tilde g}_0,z_1^{p-\delta_1}{\tilde g}_1,\cdots,z_{n-1}^{p-\delta_{n-1}}{\tilde g}_{n-1}\right)\cr
={1\over w_0^{np}}
W\left(w_0^{p-\delta_0}{\hat g}_0,w_1^{p-\delta_1}{\hat g}_1,\cdots,w_{n-1}^{p-\delta_{n-1}}{\hat g}_{n-1}\right),\cr}
$$
where
$$
{\hat g}_j(w_0,\cdots,w_{n-1})=g_j(w_0,\cdots,w_{n-1},1)
$$
for $0\leq j\leq n$.  Thus
$$
\displaylines{
{1\over\prod_{j=1}^n z_j^{p-\delta_j-n+1}}
W\left({\tilde g}_0,z_1^{p-\delta_1}{\tilde g}_1,\cdots,z_{n-1}^{p-\delta_{n-1}}{\tilde g}_{n-1}\right)\cr
=
{w_0^{p-\sum_{j=0}^n\delta_j-(n+1)(n-1)}
\over\prod_{j=1}^n w_j^{p-\delta_j-n+1}}
\left\{{1\over w_0^{p-\delta_0-n+1}}
W\left(w_0^{p-\delta_0}{\hat g}_0,w_1^{p-\delta_1}{\hat g}_1,\cdots,w_{n-1}^{p-\delta_{n-1}}{\hat g}_{n-1}\right)\right\}\cr}
$$
which is holomorphic on the whole hypersurface and vanishes on an ample
divisor, because $p>(n+1)(n-1)+\sum_{j=0}^n\delta_j$.  We conclude 
from Theorem 1 that
the Wronskian
$$
W\left(w_0^{p-\delta_0}{\hat g}_0,w_1^{p-\delta_1}{\hat g}_1,\cdots,w_{n-1}^{p-\delta_{n-1}}{\hat g}_{n-1}\right)
$$
must be identically zero on the image of $f$ (more precisely on the image
of $d^{n-1}f$) and
there is a nontrivial linear relation among
$$
x_1^{p-\delta_1}g_1(x_0,\cdots,x_n),\cdots,x_n^{p-\delta_n}g_n(x_0,\cdots,x_n)
$$
on the image of $f$.  Q.E.D.

\medskip
\noindent
{\it Proposition 3 (the Borel lemma for high powers of entire functions).}  
Let $n\geq 2$ and $p>(n-1)(n+1)$ be integers.  Let
$f_0,\cdots ,f_n$ be entire functions on  ${\bf C}$ 
satisfying  $f_0 + \cdots  + f_n \equiv  0$ such that
$f_j=g_j^p$ for some entire function $g_j$ for $1\leq j\leq n$.  Then after 
relabelling the set  $f_0,\cdots ,f_n$,  one can divide up the set  
$\{0,\cdots,n\}$  into  $q$  disjoint subsets  
$\{\ell_0,\cdots,\ell_1-1\}, 
\{\ell_1,\cdots,\ell_2-1\},\cdots,\{\ell_{q-1},\cdots,\ell_q-1\}$  
with $0=\ell_0<\ell_1<\cdots<\ell_q=n+1$ and one can find constants  
$c_{\mu,j}$ ($0\leq\mu<q$ and $\ell_\mu<j<\ell_{\mu+1}$)  
such that  $f_j\equiv c_{\mu,j}f_{\ell_\mu}$ for  
$\ell_\mu<j<\ell_{\mu+1}$ and 
$\sum^{\ell_{\mu+1}-1}_{j=\ell_\mu}f_j\equiv
\left(1+\sum^{\ell_{\mu+1}-1}_{j=\ell_\mu+1} 
c_{\mu,j}\right)f_{\ell_\mu}\equiv 0$.
  
\medskip
\noindent
{\it Proof.}  This follows from Proposition 2 by using
$\delta_j=0$ for $0\leq j\leq n$ and use induction on $n$.  

\medskip
Alternatively one can also argue directly by using Cartan's version of 
the Second Main Theorem with truncated counting function instead of
Proposition 2 as follows.

\medskip
Let  $f:{\bf C}\rightarrow{\bf P}_n$ be a 
nonconstant holomorphic map whose image is
contained in $X$.  Consider the map $\Phi:{\bf P}_n\rightarrow  
{\bf P}_{n-1}$ defined with the homogeneous coordinates  
$[(g_1)^p,\cdots,(g_n)^p]$.  Let  $H_\mu$  
($1\leq\mu\leq n-1$)  be the coordinate hyperplanes of  ${\bf P}_{n-1}$.  
Let  $H^\prime$  be the hyperplane in  ${\bf P}_{n-1}$ defined by the vanishing of 
the sum of the homogeneous coordinates of ${\bf P}_{n-1}$.  Since the image of
$f$ lies in $X$, the pullback by  $\Phi\circ f$  of the defining function 
of  $H^\prime$  is the same as the pullback by  $f$  of $-(f^*s_0)^p$.  A point  $P$
of  $H^\prime$  is assumed by  $\Phi \circ f$  at some point  $z_0$ of  ${\bf C}$  
if and only if  $-(f^*s_0)^p$ vanishes at  $z_0$ and, in that case, it must 
automatically vanish to order at least  $p$  and so that the point $P$  of  $H^\prime$  is
assumed by  $\Phi\circ f$  with multiplicity at least  $p$.  As a consequence
the truncated counting function  $N_{n-1}(r,\Phi\circ f,H^\prime)$ is no
more than ${n-1\over p}
N(r,\Phi\circ f,H^\prime)$.  Here the truncated counting function
$N_{n-1}(r,\Phi\circ f,H^\prime)$ means that multiplicities higher than
$n-1$ are replaced by multiplicites equal to $n-1$.
The same argument holds for  $H_j$ ($1\leq j\leq N$)  
instead of  $H^\prime$.  Unless the image of  $\Phi\circ f$  is contained in a 
hyperplane of  ${\bf P}_{n-1}$,  we know from Cartan's Second Main Theorem with 
truncated counting function that
$$
\displaylines{
T(r,\Phi\circ f)\leq N_{n-1}(r,\Phi\circ f,H^\prime)+\sum^n_{j=1} 
N_{n-1}(r,\Phi\circ f,H_j)+O(\log T(r,\Phi\circ f))
\cr
\leq{n-1\over p}(N(r,\Phi \circ f,H^\prime)+\sum ^n_{j=1} 
N(r,\Phi\circ f,H_j))+O(\log T(r,\Phi\circ f))
\cr
\leq{n-1\over p}(n+1)T(r,\Phi\circ f)+O(\log T(r,\Phi\circ f)),
\cr
}
$$
which is a contradiction if $p>(n+1)(n-1)$.  So we conclude that the image 
of $f$ is contained in the zero-set of $\sum^n_{j=1}\lambda_j(g_j)^p=0$  
for some  $\lambda_j\in{\bf C}$ ($1\leq j\leq n$) not all zero.  
Now we use induction on $n$. Q.E.D.

\medskip
Now we introduce the argument of counting dimensions of certain subvarieties
of Grassmannians.  Let $2\leq k\leq m-1$ and $\ell\geq 2$.
Let ${\cal G}$ be the Grassmannian of ${\bf C}^k$ in ${\bf C}^m$.  For
hyperplanes $H_1,\cdots,H_\ell$ in ${\bf C}^m$ let
${\cal G}_{H_1,\cdots,H_\ell}$ be the subvariety of
$\cal G$ consisting of all
$W\in{\cal G}$ such that the dimension over ${\bf C}$ of the 
linear space generated by $H_1|W,\cdots,H_\ell|W$ is no more than $1$.
For generic hyperplanes $H_1,\cdots,H_\ell$ in ${\bf C}^m$ the dimension
over ${\bf C}$ of ${\cal G}_{H_1,\cdots,H_\ell}$ is equal to the dimension
of the Grassmann ${\cal G}^\prime$
of all ${\bf C}^k$ in ${\bf C}^{m-\ell+1}$ plus $\ell-1$.
The reason is as follows.  For $W\in{\cal G}_{H_1,\cdots,H_\ell}$,
if the linear space generated by $H_1|W,\cdots,H_\ell|W$ 
is precisely of dimension $1$, after relabelling $H_1,\cdots,H_\ell$, 
we have constants
$c_2,\cdots,c_\ell$ such that $H_j=c_jH_1$ for $2\leq j\leq\ell$ and
$W$ is an element of
the Grassmannian ${\cal G}^\prime$ of all ${\bf C}^k$ in 
$$
{\bf C}^{m-\ell+1}={\bf C}^m\cap\{H_2=c_2H_1,\cdots,H_\ell=c_\ell H_1\},
$$
and the freedom of choices for $c_2,\cdots,c_\ell$ gives the additional
$\ell-1$ dimensions.
If all $H_1|W,\cdots,H_\ell|W$ are identically zero, then $W$ is an element of
the Grassmannian ${\cal G}^{\prime\prime}$ of all ${\bf C}^k$ in 
$$
{\bf C}^{m-\ell}={\bf C}^m\cap\{H_1=\cdots=H_\ell=0\}
$$ and ${\cal G}^{\prime\prime}\subset{\cal G}^\prime$.
Thus for generic hyperplanes $H_1,\cdots,H_\ell$ in ${\bf C}^m$ the 
codimension of ${\cal G}_{H_1,\cdots,H_\ell}$ in ${\cal G}$ is
$(k-1)(\ell-1)$.

\medskip
Let $q\geq 1$ and $\ell_\nu\geq 2$ for $1\leq\nu\leq q$.
For generic hyperplanes
$$
H^{(1)}_1,\cdots,H^{(1)}_{\ell_1},
\cdots,
H^{(q)}_1,\cdots,H^{(q)}_{\ell_q},
$$
the codimension over ${\bf C}$ of
$$
{\cal G}_{H^{(1)}_1,\cdots,H^{(1)}_{\ell_1}}
\cap\cdots\cap
{\cal G}_{H^{(q)}_1,\cdots,H^{(q)}_{\ell_q}}
$$
is equal to
$(k-1)\sum_{\nu=1}^q(\ell_\nu-1)$ when
$(k-1)\sum_{\nu=1}^q(\ell_\nu-1)\leq k(m-k)$ and is empty
$(k-1)\sum_{\nu=1}^q(\ell_\nu-1)>k(m-k)$
Let $N=\sum_{\nu=1}^q \ell_\nu$.

\medskip
\noindent
{\it Proposition 4}.  For $N\geq 4m-7$ and for generic hyperplanes 
$$
H^{(1)}_1,\cdots,H^{(1)}_{\ell_1},
\cdots,
H^{(q)}_1,\cdots,H^{(q)}_{\ell_q},
$$
in ${\bf C}^m$ the subvariety
$$
{\cal G}_{H^{(1)}_1,\cdots,H^{(1)}_{\ell_1}}
\cap\cdots\cap
{\cal G}_{H^{(q)}_1,\cdots,H^{(q)}_{\ell_q}}
$$
of $\cal G$ is empty.

\medskip
\noindent
{\it Proof}.  It follows from $\ell_\nu\geq 2$ that
$N=\sum_{\nu=1}^q \ell_\nu\geq 2q$.  The codimension of
the subvariety
$$
{\cal G}_{H^{(1)}_1,\cdots,H^{(1)}_{\ell_1}}
\cap\cdots\cap
{\cal G}_{H^{(q)}_1,\cdots,H^{(q)}_{\ell_q}}
$$
in $\cal G$ is 
$$
(k-1)\sum_{\nu=1}^q(\ell_\nu-1)=(k-1)(N-q)\geq
(k-1){N\over 2}>k(m-k)=\dim_{\bf C}{\cal G}.
$$
Q.E.D.

\medskip
\noindent
{\it Proof of Theorem 3}. Suppose the contrary and we are going to derive a contradiction.
Let $f:{\bf C}\rightarrow{\bf P}_n$ be a nonconstant holomorphic map whose
image lies in the hypersurface $\sum_{j=1}^N H_j^p=0$.  Without loss of
generality we can assume that the image of $f$ does not lie in the
the hyperplane $H_j=0$ for any $1\leq j\leq N$, otherwise we argue instead
with the hypersurface $\sum_{j\in J}H_j^p=0$, where $J$ is the set of all
$1\leq j\leq N$ with the property that the image of $f$ does not lie in the
the hyperplane $H_j=0$.  Let
$\tilde f:{\bf C}\rightarrow{\bf C}^{n+1}$ be a lifting of $f$.  Let
$W$ be the linear subspace of ${\bf C}^{n+1}$ which is the linear span
of the image of $\tilde f$.  Let $k$ be the complex dimension of
$W$.  Since $f$ is nonconstant, we know that $k$ is at least $2$.  
Let $m=n+1$. By the
Borel lemma for high powers of entire functions, we have a partition
$$
\{H^{(1)}_1,\cdots,H^{(1)}_{\ell_1}\}
\cup
\cdots
\cup
\{H^{(q)}_1,\cdots,H^{(q)}_{\ell_q}\}
$$
of $\{H_0,\cdots,H_N\}$ with $\ell_j\geq 2$ for $1\leq j\leq q$
so that the complex dimension of the linear space spanned by
$H^{(\nu)}_1|W,\cdots,H^{(\nu)}_{\ell_\nu}|W$ is at most $1$ for
$1\leq\nu\leq q$, contradicting the
preceding lemma which states that
the subvariety
$$
{\cal G}_{H^{(1)}_1,\cdots,H^{(1)}_{\ell_1}}
\cap\cdots\cap
{\cal G}_{H^{(q)}_1,\cdots,H^{(q)}_{\ell_q}}
$$
of $\cal G$ is empty.  Q.E.D.

\bigskip
\noindent
{\bf \S5. Hyperbolic Surface of Degree 11}

\medskip
We now prove Theorem 4.  Assume $n\geq 11$ and denote by $S$ the surface
defined by the equation
$$
x_0^n+x_1^n+x_2^n+x_3^{n-2}g(x_0,x_1,x_2,x_3)=0.
$$
Suppose we have a nonconstant
holomorphic map $f:{\bf C}\rightarrow S$ and we are going to derive
a contradiction.
By Proposition 2, we have a nontrivial linear relation among
$x_0^n,x_1^n,x_2^n$ on the image of ${\bf C}$, 
which we can assume without loss of generality
to be $x_0^n=c_1x_1^n+c_2x_2^n$.  
When both $c_1,c_2$ are nonzero, the image of ${\bf C}$ lies in the
Fermat curve $x_0^n=c_1x_1^n+c_2x_2^n$ and we end up with
$[x_0,x_1,x_2]=\hbox{constant}$ on the image of $\bf C$.  We cannot have
all three $x_0,x_1,x_2$ identically zero on the image of $\bf C$, 
because the assumption $g(0,0,0,x_3)=x_3^2$
would imply that $x_3$ is identically zero on the image of $\bf C$
as well.  Suppose without loss of generality that $x_0$ is not identically 
zero.  Then we end up with
$$
b_0x_0^n+x_3^{n-2}(x_3^2+b_1x_3x_0+b_2x_0^2)
$$
for some constants $b_0,b_1,b_2$, implying that the image of $f$ is constant.

\medskip
Before we continue further, we make the trivial observation that,
for a quadratic polynomial $h(y)=Ay^2+By+C$ of a single variable $y$,
$$
{(h^\prime)^2\over 2h^{\prime\prime}}-h={B^2-4AC\over 4A},
$$
where $h^\prime$ and $h^{\prime\prime}$ denote respectively the first
and second derivatives of $h$ with respect to $y$.

\medskip
Now assume that $c_2=0$ and $x_0^n=c_1x_1^n$ on the nonconstant image of
${\bf C}$, we conclude that $x_0=c_1^{1\over n}x_1$ for some $n^{\hbox{th}}$
root of $c_1$ on the image of ${\bf C}$.  Thus the curve $C$ defined by
$$
\displaylines{
x_0=c_1^{1\over n}x_1,\cr
(1+c_1)x_1^n+x_2^n+x_3^{n-2}g(c_1^{1\over n}x_1,x_1,x_2,x_3)=0\cr
}
$$
contains the image of ${\bf C}$.  Let $U_3$ be the affine open subset
$\{x_3\not=0\}$.  Then $C\cap U_3$ is defined, in terms of the affine 
coordinates $\zeta_j={x_j\over x_3}$ ($0\leq j\leq 2$) by
$$
\displaylines{
\zeta_0=c_1^{1\over n}\zeta_1,\cr
(1+c_1)\zeta_1^n+\zeta_2^n+
g(c_1^{1\over n}\zeta_1,\zeta_1,\zeta_2,1)=0.\cr
}
$$
We distinguish between two cases.  First we consider the
case $1+c_1=0$.  In that case we have
$$
\displaylines{
\zeta_0=(-1)^{1\over n}\zeta_1,\cr
\zeta_2^n+g((-1)^{1\over n}\zeta_1,\zeta_1,\zeta_2,1)=0.\cr
}
$$
Denote $g((-1)^{1\over n}\zeta_1,\zeta_1,\zeta_2,1)$ by 
$h(\zeta_1,\zeta_2)$ and let
$$
h(\zeta_1,\zeta_2)=gA(\zeta_2)\zeta_1^2+B(\zeta_2)\zeta_1+C(\zeta_2).
$$
Then
$$
\displaylines{
0=\zeta_2^n+h(\zeta_1,\zeta_2)\cr
=\zeta_2^n+A(\zeta_2)\left(\zeta_1+{B(\zeta_2)\over 4A(\zeta_2)}\right)^2
-{B(\zeta_2)^2-4A(\zeta_2)C(\zeta_2)\over 4A(\zeta_2)}\cr
=\zeta_2^n+A(\zeta_2)\left(\zeta_1+{B(\zeta_2)\over 4A(\zeta_2)}\right)
-\left(
{\left({\partial h\over\partial\zeta_1}\right)^2\over
2{\partial^2h\over\partial\zeta_1^2}}-h
\right)\cr
}
$$
and
$$
A(\zeta_2)\left(\zeta_1+{B(\zeta_2)\over 4A(\zeta_2)}\right)^2
=-\zeta_2^n+\left(
{\left({\partial h\over\partial\zeta_1}\right)^2\over
2{\partial^2h\over\partial\zeta_1^2}}-h
\right).
$$
Since the polynomial 
$$
-\zeta_2^n+\left(
{\left({\partial h\over\partial\zeta_1}\right)^2\over
2{\partial^2h\over\partial\zeta_1^2}}-h
\right)
$$
in $\zeta_2$ has $n$ distinct roots by assumption,
it follows that the normalization of the hyperelliptic Riemann surface defined by
$$
\zeta_2^n+g((-1)^{1\over n}\zeta_1,\zeta_1,\zeta_2,1)=0
$$
has genus equal to $\lceil{n\over 2}\rceil-1$, where $\lceil\cdot\rceil$
denotes the round-up.  This contradicts the nonconstancy of $f$.

\medskip
Now we consider the case $1+c_1\not=0$.  Then from the 
second equation in
$$
\displaylines{
x_0=c_1^{1\over n}x_1,\cr
(1+c_1)x_1^n+x_2^n+x_3^{n-2}
g(c_1^{1\over n}\zeta_1,\zeta_1,\zeta_2,1)=0\cr
}
$$
we conclude that $x_1^n,x_2^n$ have a nontrivial linear relation on the
image of ${\bf C}$.  Again we conclude that the image of ${\bf C}$ is constant.
This concludes the proof of Theorem 4.

\medskip
\noindent
{\it Remark.} In Theorem 4, one can also easily formulate a necessary 
and sufficient condition on $g(x_0,x_1,x_2,x_3)$ so that 
the surface is hyperbolic for $n\geq 11$.

\bigskip
\noindent
{\bf References}

\medskip
\noindent
[A41] L. Ahlfors, The Theory of meromorphic curves, Acta Soc. Sci. Fenn. Nova 
Ser. A {\bf 3}(4) (1941), 171-183.

\medskip
\noindent
[AS80] K. Azukawa and M. Suzuki, Some examples of algebraic degeneracy and 
hyperbolic manifolds, Rocky Mountain J. Math. {\bf 10} (1980), 655-659.

\medskip
\noindent
[B26] A. Bloch, Sur les syst\`emes de fonctions uniformes satis\-faisant \`a 
l'\'equat\-ion d'une vari\'et\'e alg\'ebrique dont l'irr\'egularit\'e 
d\'epasse la dimension.  J. de Math. {\bf 5} (1926), 19-66.

\medskip
\noindent
[BG77] R. Brody and M. Green, A family of smooth hyperbolic hypersurfaces of  
${\bf P}_3$, Duke Math. J. {\bf 44} (1977), 873-874.

\medskip
\noindent
[D95] J.-P. Demailly, Lecture in Nancy 1995 and a series of lectures at
Santa Cruz 1995.

\medskip
\noindent
[DL96] G. Dethloff and S. Lu, Lecture in MSRI, 1996.

\medskip
\noindent
[E96] J. El Goul, Algebraic Families of Smooth Hyperbolic surfaces of 
Low Degree in $P^3_C$, Manuscripta Math. {\bf 90} (1996), 521-532.

\medskip
\noindent
[GG79] M. Green and P. Griffiths, Two applications of algebraic geometry to 
entire holomorphic mappings, The Chern Symposium 1979, Proc. Internat. Sympos.,
Berkeley, 1979, Springer-Verlag 1980.

\medskip
\noindent
[K80] Y. Kawamata, On Bloch's conjecture, Invent. Math. {\bf 57} (1980), 97-100.

\medskip
\noindent
[MN94] K. Masuda and J. Noguchi, A construction of hyperbolic hypersurfaces of
${\bf P}^n({\bf C})$, Preprint 1994.

\medskip
\noindent
[McQ96] M. McQuillan, A new proof of the Bloch conjecture, J. Alg. Geom. 
{\bf 5} (1996), 107-117.

\medskip
\noindent
[Na89] A. Nadel, Hyperbolic surfaces in  ${\bf P}^3$, Duke Math. J. {\bf 58} (1989), 
749-771.

\medskip
\noindent
[NO90] J. Noguchi and T. Ochiai, {\it Geometric Function Theory in Several Complex 
Variables}, Transl. Math. Mon. {\bf 80}, Amer. Math. Soc., Providence, R.I. 1990.

\medskip
\noindent
[Oc77] T. Ochiai, On holomorphic curves in algebraic varieties with ample 
irregularity, Invent. Math. {\bf 443} (1977), 83-96.

\medskip
\noindent
[S87] Y.-T. Siu, Defect relations for holomorphic maps between
spaces of different dimensions, Duke Math.\ J.\ {\bf 55} (1987),
213--251.

\medskip
\noindent
[S95] Y.-T. Siu, Hyperbolicity problems in function theory, in ``Five Decades as 
a Mathematician and Educator - on the 80th birthday of Professor
Yung-Chow Wong'' ed. Kai-Yuen Chan and Ming-Chit Liu, World Scientific:
Singapore, New Jersey, London, Hong Kong, 1995, pp.409-514.

\medskip
\noindent
[SY96a] Y.-T. Siu and S.-K. Yeung, Hyperbolicity of the complement of 
a generic smooth curve of high degree 
in the complex projective plane, Inventiones Math. {\bf 124} (1996), 
573-618.

\medskip
\noindent
[SY96b] Y.-T. Siu and S.-K. Yeung, A generalized Bloch's theorem and the 
hyperbolicity of the complement of an ample divisor in an abelian variety,
Math. Ann., to appear.

\medskip
\noindent
[W80] P.-M. Wong, Holomorphic mappings into Abelian varieties, Amer. J. Math. 
{\bf 102} (1980), 493-501.

\medskip
\noindent
[Z89] M. Zaidenberg, Stability of hyperbolic embeddedness and construction of 
examples, Math. USSR Sbornik {\bf 63} (1989), 351-361.

\medskip
\noindent
Authors' addresses:

\medskip
\noindent
Yum-Tong Siu, Department of Mathematics, Harvard University, Cambridge, MA 02193
(e-mail: siu@math.harvard.edu).

\medskip
\noindent
Sai-Kee Yeung, Department of Mathematics, Purdue Univer\-sity, West Lafay\-ette, IN 47907
(e-mail yeung@math.purdue.edu).

\end{document}